\documentclass[11pt]{article}
\usepackage[final]{epsfig}
\usepackage{graphics}
\usepackage{amsmath}
\usepackage{amsfonts}
\usepackage{latexsym}
\usepackage{epsfig}
\usepackage{amssymb}
\usepackage{graphicx}
\usepackage{url}
\usepackage{epstopdf}
\usepackage[utf8]{inputenc}
\usepackage[russian,english]{babel}

\newtheorem{lemma}{Lemma}[section]

\newtheorem{remark}[lemma]{Remark}

\newtheorem{theorem}{Theorem}

\begin{document}
\newcommand{\eps}{{\varepsilon}}
\newcommand{\proofend}{$\Box$\bigskip}
\newcommand{\CC}{{\mathbb C}}
\newcommand{\Q}{{\mathbb Q}}
\newcommand{\R}{{\mathbb R}}
\newcommand{\Z}{{\mathbb Z}}
\newcommand{\RP}{{\mathbb {RP}}}
\newcommand{\CP}{{\mathbb {CP}}}
\newcommand{\Tr}{\rm Tr}
\def\proof{\paragraph{Proof.}}

\title{Centers of mass of Poncelet polygons, 200 years after}

\author{Richard Schwartz\footnote{
Department of Mathematics, Brown University, Providence, RI 02912;
res@math.brown.edu} \and 
Serge Tabachnikov\footnote{
Department of Mathematics,
Penn State University,
University Park, PA 16802;
tabachni@math.psu.edu}
}

\date{}
\maketitle

\section{A  letter from Saratov} \label{story}

During his last trip to Moscow, the second author of this article came into possession of a remarkable mathematical letter. The custodian of the letter, a Russian businessman X who wished to remain anonymous,
presented the letter to Tabachnikov at the end of his lecture on configuration theorems in projective geometry \cite{ST} (like many successful contemporary Russian entrepreneurs, Mr. X has a  degree in mathematics). 

Mr. X explained that the letter had been written by his great-great-great grandfather, Konstantin Shestakov, shortly after he had been discharged from the Russian Army during the Napoleonic Wars. It appears that K. Shestakov had befriended the famous French mathematician, Jean-Victor Poncelet, then a prisoner of war in the Russian town of Saratov, and had been drawn into geometry by him.\footnote{According to the family history, Shestakov was well educated: he attended Saint Petersburg Engineering School.} The letter, dated  fall of 1814 and apparently never mailed,  was addressed to Shestakov's younger brother, Alexander, who had attended the Kazan Gymnasium with N. I. Lobachevsky.

We reproduce the first paragraph of the letter in the original Russian (in the modern orthography) and give its full translation  into English.

\begin{quote}
{\small \it
\begin{otherlanguage}{russian}
Милый Саша!

Ты конечно слышал, что раны мои, полученные при битве под Смоленском, вынудили меня выйти в отставку, и что я постепенно восстанавливаю силы свои в Саратове. Это сонный город на Волге, где окромя роскошных обедов да балов у здешнего губернатора г. Панчулидзева и делать вовсе нечего. Служба моя состоит и в надзоре за здешней тюрьмой, однако заключённые ведут себя смирно и не нуждаются в постоянном присмотре. Пишу сообщить о своём житье-бытье, а также и с просьбой переслать математическую часть сего письма однокашнику твоему по гимназии Николаю Лобачевскому; как слыхал я, он теперь адъюнкт-профессор в Казани... 
 \end{otherlanguage}
}
\end{quote}

\begin{quote}
{\small \it
Dearest Sasha:

You doubtless have heard that the wounds sustained at a battle near Smolensk
forced me to withdraw from the army  and that I have been recovering ever since in Saratov. Saratov is a dreary town on the Volga river which, aside from the several excellent dinner parties and balls thrown by 
Governor Panchuleedzev, affords us almost nothing interesting to do. I have some hand in overseeing the prison, but the  prisoners are peaceful and rarely need attention. I am writing this partly to tell you how I have been keeping myself, but also to ask you the favor of passing on the mathematical part of this letter to your former classmate Nikolai Lobachevsky who, as I understand it, is now an adjunct professor at the Kazan University.

I spent much of the winter holed up in my room, trying to avoid the oppressive regimen of drinking and other kinds  of enforced merriment. As you know, I prefer a life of quiet contemplation. I had never quite become used to the cruel interruption of my studies brought about by the war. I confess that my bitterness over the situation had been driving me towards madness, but then everything changed. 
I write to you now with great excitement, as a man given a divine gift.

My gift came to me in the form a prisoner, Lieutenant of the French Army Monsieur Poncelet. As luck would have it, Poncelet was assigned to me as an assistant of sorts, though the pedestrian kind of work required of him hardly touched on his enormous talents. It became clear almost immediately that he was nothing at all like our other prisoners. Quiet and unassuming, with a temperament much like my own, Poncelet preferred to keep his own company. The relatively undemanding work he did for me afforded him to pursue his true passion, which was geometry.

Poncelet explained to me that he was whiling away the time by reconstructing from first principles the theorems of geometry he had learned during his student days at   \'Ecole Polytechnique in Paris. He said more than once that the great beauty of mathematics is that it can be created whole from very few principles. I knew something of geometry myself and I gradually became engaged in the project of this brilliant  man. In the space of several months I went from casual observer to eager student to active participant.

You have to understand that Poncelet was no ordinary student! As time passed, I could see that he was not just trying to recreate the lessons of his masters, as I had first thought, but rather that he was striking out into the unknown. My excitement was almost as great as his when he explained to me his astonishing discovery about conic sections.

This is the part of the letter I request that you pass on to Lobachevsky. 
I recall that your mathematics teacher, the venerable Grigorii Ivanovich, always praised Nikolai for his aptitude for mathematics. I trust that you will be able to follow this but, in any case, Lobachevsky certainly will, and I hope that it will interest him.
 
Suppose that $A$ and $B$ are elliptical conic sections with $A$ surrounding $B$. Imagine that one has a point $P_1$ on $A$, then draws a segment $P_1P_2$ which is tangent to $B$ and connects $P_1$ to $P_2$, also on $A$. Now repeat the construction, drawing the segment $P_2 P_3$ tangent to $B$ so that $P_3$ is on $A$. Imagine that this construction is repeated some number of times, say 100, so that $P_{101}$ is the same point as $P_1$. One has produced a kind of closed polygonal figure involving points $P_1, P_2$, etc., which ends up where it starts after 100 steps.

Poncelet's great discovery is that the same construction, starting with a
different point on $A$, will also repeat after 100 steps. That is, if one moves
$P_1$ to the new point $P'_1$ on $A$, and then produces points $P'_2, P'_3$, etc. then 
it will again happen that $P'_{101} = P'_1$! You have to understand, finally, that 
I have only used 100 as an example. Whether the figure repeats itself or not, and after how many steps, depends on the choice of $A$ and $B$. Poncelet encouraged me to think about this construction, as he did, as a kind of spinning polygon. He imagined moving the point $P_1$ continuously, so that the other points would move as well. Poncelet's great discovery is that this configuration of (say) 100 points and 100 segments spins around the conics $A$ and $B$ and remains intact.

Poncelet left the camp in June and I do not know what became of him, but I continued his work. I could not get the idea of the spinning configuration out of my mind and eventually I thought to ask some questions of my own. For simplicity I considered the case when $A$ and $B$ were both circles, with $A$ surrounding $B$ and $B$ being set in a position that was somewhat offset from the center of $A$. By adjusting the center and position of $B$ I could control the particulars of Poncelet's construction. After some amount of fussing around, I settled on a location and size of $B$ which caused the configuration to close up after 5 steps, making a star-like pattern.

I became interested in the question of how the center of gravity of the figure changed as it spun around. You might say that I approached this in an experimental way, drawing the figure on the page and then approximating the center of gravity by a kind of trick of drawing the points together. I wish that I could include diagrams of this laborious method, but perhaps they would mean nothing to you. Suffice it to say that my calculations seemed to show a promising result but were not accurate enough to convince me.

You probably remember my passion for precise reasoning! Eventually I joined many sheets together, making a kind of poster on which to do my calculations. I remember spending the bulk of a week scribbling on this poster, night after night, as gradually the beautiful answer revealed itself. It is yet another circle! As the great Poncelet figure spins around, its center of gravity traces out a circle as well! Understand that there is nothing special about the star configuration that I settled upon; I expect that the principle should be completely general. I am convinced that this is a first-rate discovery to rival Poncelet's own.

I regret that I have so far not hit upon a proof of my discovery, in the same way that Poncelet was able to find a proof for his. Alas, perhaps if my friend were still here we could find it together. Nonetheless, I have calculated things out to such a precision that I am confident in proclaiming this as a theorem. 
I would be grateful if Lobachevsky could look into the matter.

Your brother Konstantin.
}
\end{quote}

\begin{figure}[ht]
\centering
\includegraphics[height=2in]{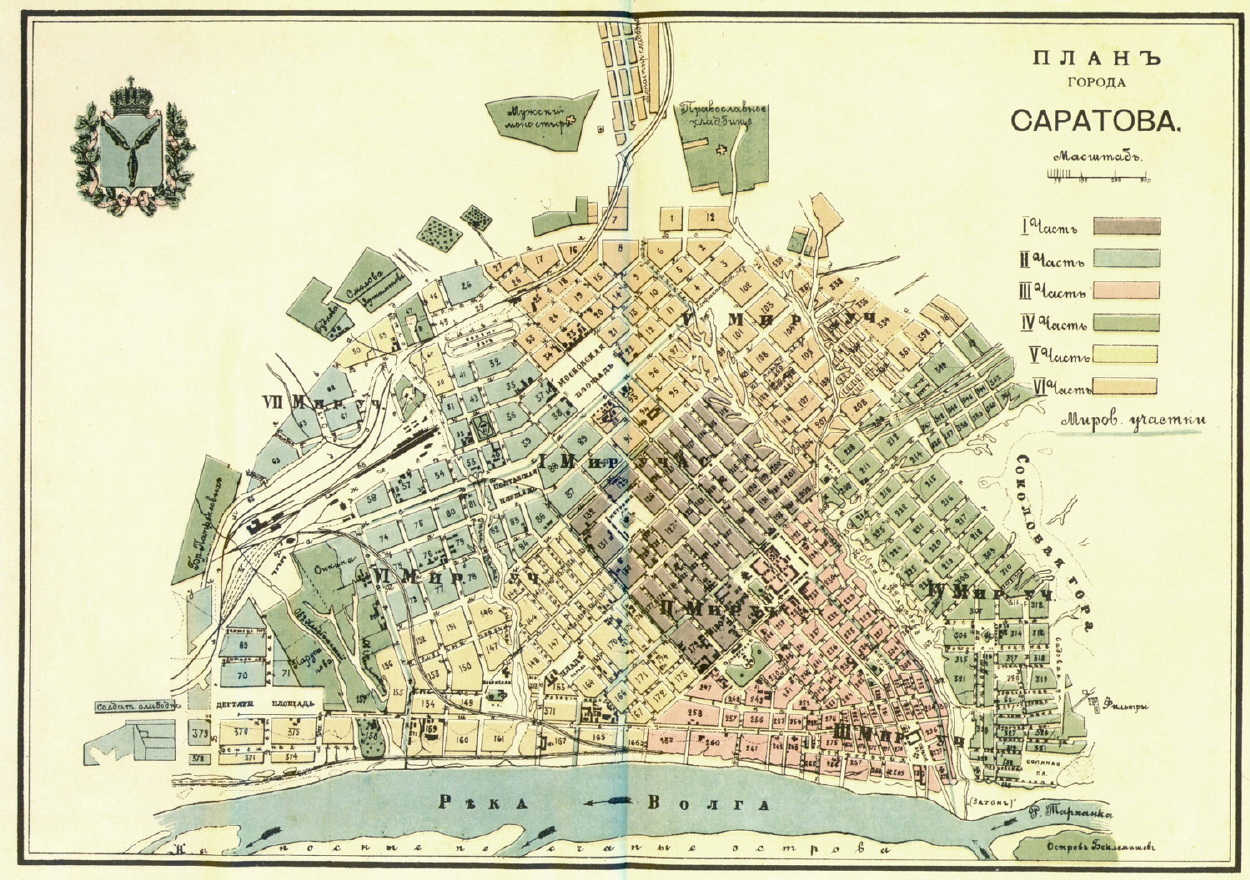}\qquad
\includegraphics[height=2in]{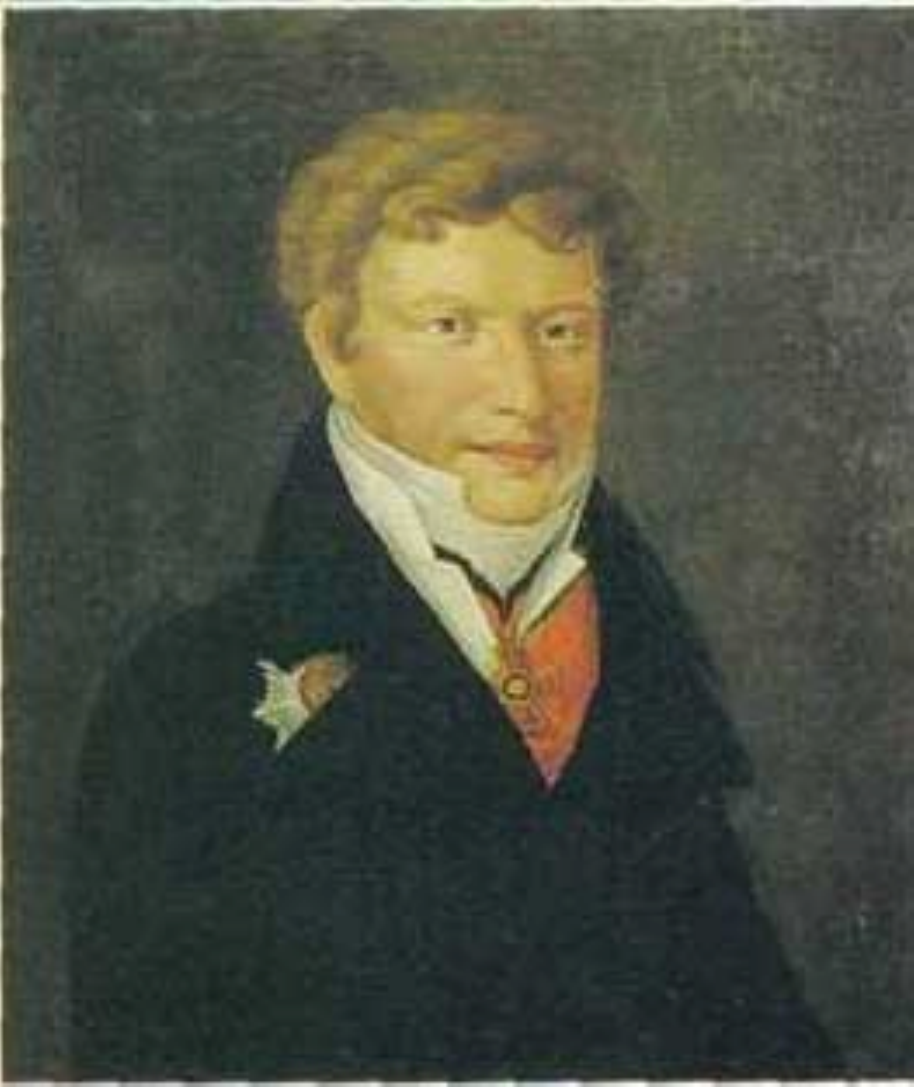}
\caption{An old map of Saratov and its Governor A. D. Panchuleedzev (1758--1834).}
\label{Saratov}
\end{figure}

To put this letter into a historical context, Jean-Victor Poncelet participated in Napoleon's invasion of Russia and was wounded at the Battle of Krasnoi during Napoleon's retreat from Moscow in November of 1812. He 
 spent more than a year in Russian captivity in Saratov. There he developed projective geometry. Later, his {\it Saratov notebooks} became part of his {\it Trait\'e des propri\'et\'es projectives des figures} (1822)
and   {\it Applications d'analyse et de g\'eom\'etrie} (1862). In particular, Poncelet discovered his celebrated porism  during the Saratov exile.

N. I. Lobachevsky grew up in Kazan where he attended gymnasium\footnote {His teacher of mathematics was  G. I. Kartashevsky.} and then the newly established Kazan University. The whole career of Lobachevsky was spent at this university where he served from 1814, as an adjunct professor, until 1846 (elected Chancellor in 1827). 

\begin{figure}[ht]
\centering
\includegraphics[height=2in]{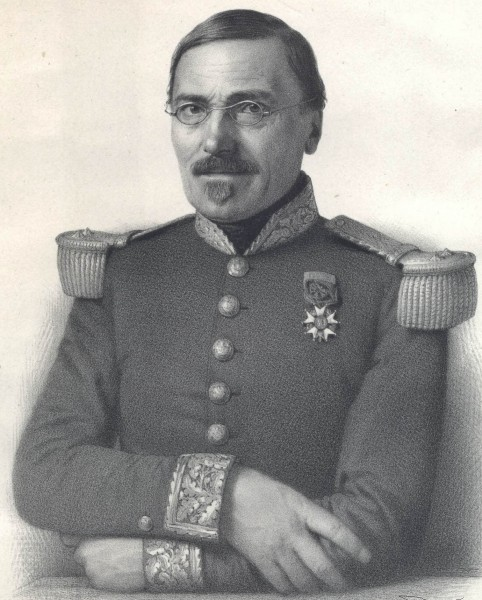}\qquad\qquad
\includegraphics[height=2in]{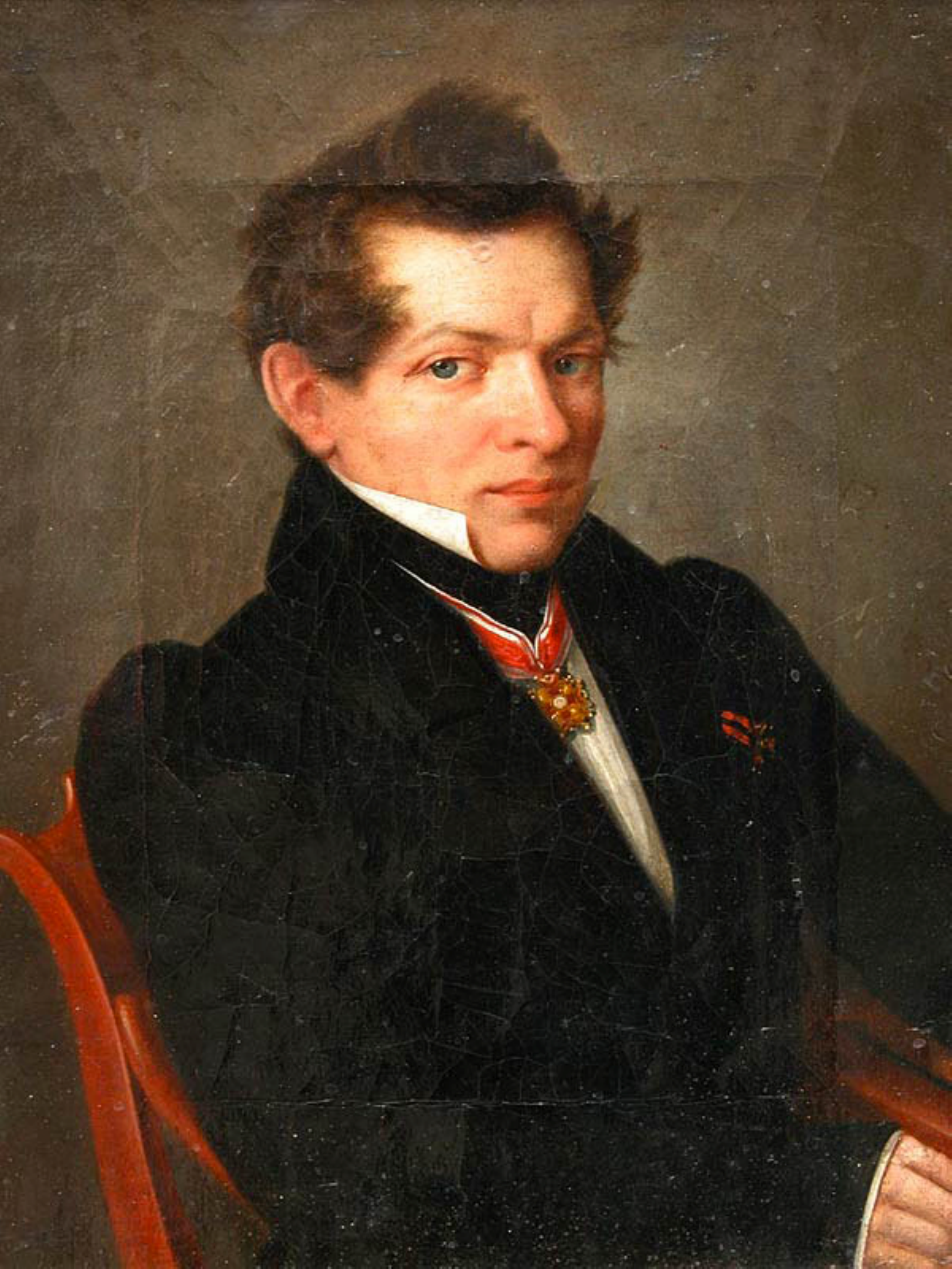}
\caption{J.-V. Poncelet (1788--1867) and N. I. Lobachevsky (1792--1856).}
\label{mathematicians}
\end{figure}

The literature on the Poncelet theorem and its ramifications is vast; see, e.g., \cite{BB,BKOR,Do,DR1, DR2,Fl} and the references therein. Still, after having mulled over
Shestakov's letter at some length, we feel that his result appears to be a
new and surprising addition to the Poncelet theorem.  What makes the result
especially surprising is that it makes a statement about the center of mass -- an
affine geometry construction -- in the context of projective geometry.
The one result we know in this direction is Weill's Theorem, which
deals with the centers of mass of Poncelet polygons in a special case.
One could view Shestakov's Theorem as a generalization of Weill's Theorem,
though Shestakov's result does not quite imply Weill's result.

In what follows, we will present and sketch modern proofs for
two versions of Shestakov's result.  
For reasons we will explain, we think that he
most likely had the first statement in mind, but it seems reasonable to present 
the second statement as well, which is a variant.  At the end, and for the sake
of completeness, we will sketch a proof of Weill's Theorem.

\section{The theorem} \label{mainform}

\begin{figure}[ht]
\centering
\includegraphics[height=2.3in]{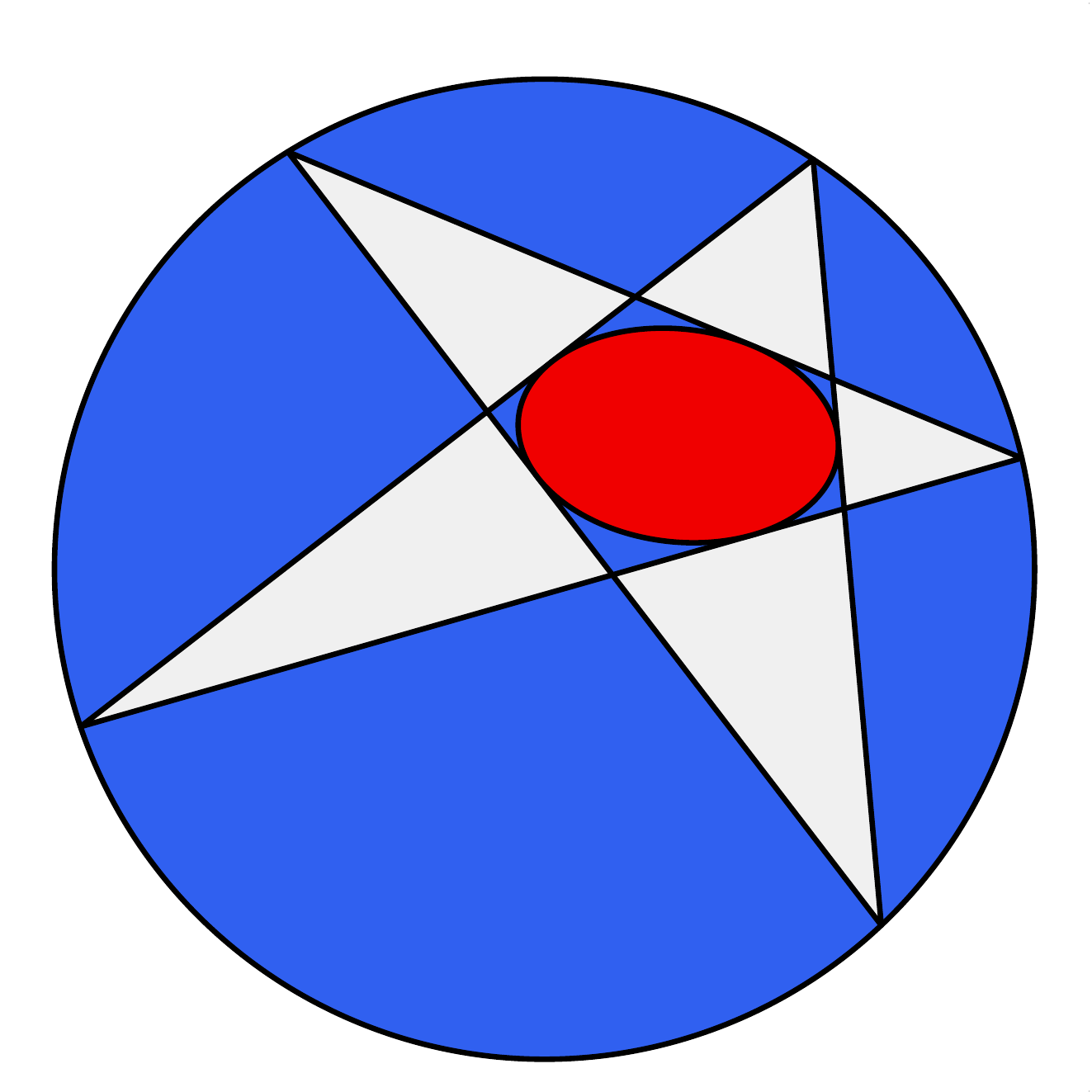}
\caption{A Poncelet pentagon.}
\label{pentagon}
\end{figure}

Figure \ref{pentagon} features an instance of the Poncelet porism -- exactly
the kind Shestakov said that he had ``fussed around to arrange''.
The Poncelet theorem asserts that every point 
of the outer conic is a vertex of such a polygon.
The interested reader can use\footnote{You can download the
program at http://www.math.brown.edu/$\sim$res/Java/PORISM.tar}
 our program and see
the kind of animation envisioned by Poncelet.

Let $P=(A_1\ldots A_n)$ be a Poncelet $n$-gon, 
and let $(x_i,y_i)$ be the Cartesian coordinates 
of the vertex $A_i$, $i=1,\ldots,n$.  
There are really $3$ natural interpretations for the center of mass of
$P$:
\begin{itemize}
\item The center of mass of the vertices.  We denote this by $CM_0(P)$.
\item The center of mass of the edges, when the edges are given a
uniform density.  We denote this by $CM_1(P)$.
\item The center of mass of the $P$ when $P$ is considered as
a ``homogeneous lamina''. 
We denote this by $CM_2(P)$. When $P$ is not
embedded, $CM_2(P)$ requires the interpretation
we give in Lemma \ref{CMcoord} below.
\end{itemize}

We provide formulas for the three kinds of centers below.
We will present version of Weill's Theorem for $CM_0$ and $CM_2$.
The result fails for $CM_1$ for reasons we will mention briefly below.
We think that Shestakov most likely had $CM_0$ in mind from his
description of ``drawing the points together''.  This is the
sort of thing one would do when computing $CM_0(P)$.  

The coordinates of  these centroids are given in the next lemma. Let 
$$
\ell_i=\sqrt{(x_{i+1}-x_i)^2 + (y_{i+1}-y_i)^2}
$$
be the length of $i$th side of the polygon, and let
$d_i = (x_i y_{i+1} - x_{i+1} y_i). $ Then $L(P)=\sum_i \ell_i$ is the perimeter  of $P$, and
$A(P)=(1/2) \sum_i d_i$ is  its area, counted with sign and multiplicity.

\begin{lemma}
 \label{CMcoord}
One has:
\begin{itemize}
\item $CM_0(P) =\frac{1}{n} \sum_{i=1}^n (x_i, y_i),$
\item $CM_1(P) = \frac{1}{2L(P)} \sum_{i=1}^n \ell_i  (x_i + x_{i+1},y_i + y_{i+1}),$
\item $CM_2(P)= \frac{1}{6 A(P)} \sum_{i=1}^n d_i (x_{i}+x_{i+1},y_{i}+y_{i+1}).$
\end{itemize}
\end{lemma}

\proof Consider the triangle $O A_i A_{i+1}$ where $O$ is the origin. The centroid of this triangle is at point $(A_i + A_{i+1})/3$ and it has the area $d_i/2$. Summing up over $i$ and dividing by the total area of $P$ yields the result for $CM_2(P)$. Likewise for $CM_0(P)$ and $CM_1(P)$.
\proofend 

Here is (our interpretation of) Shestakov's Theorem.

\begin{theorem} \label{main}
Let $\gamma \subset \Gamma$ be a pair of nested ellipses that admit a 1-parameter family of Poncelet $n$-gons $P_t$.  Then both loci $CM_0(P_t)$ and $CM_2(P_t)$ are ellipses
homothetic to $\Gamma$ (or single points).
\end{theorem} 

As we mentioned above, Shestakov most likely had $CM_0(P)$ in mind.
Also, he has normalized so that the outer ellipse is a circle; in this
case the loci $CM_0(P_t)$ and $CM_2(P_t)$ are circles.  For the proof
we will also normalize this way.  Our proof relies on the fact that
$CM_0(P)$ and $CM_2(P)$ are rational expressions in the coordinates
of $P$.  This is not true for $CM_1(P)$ and, as we mentioned above,
the result fails for $CM_1(P)$.  The locus in this case is not
generally a conic section.

\begin{figure}[ht]
\centering
\includegraphics[height=2.3in]{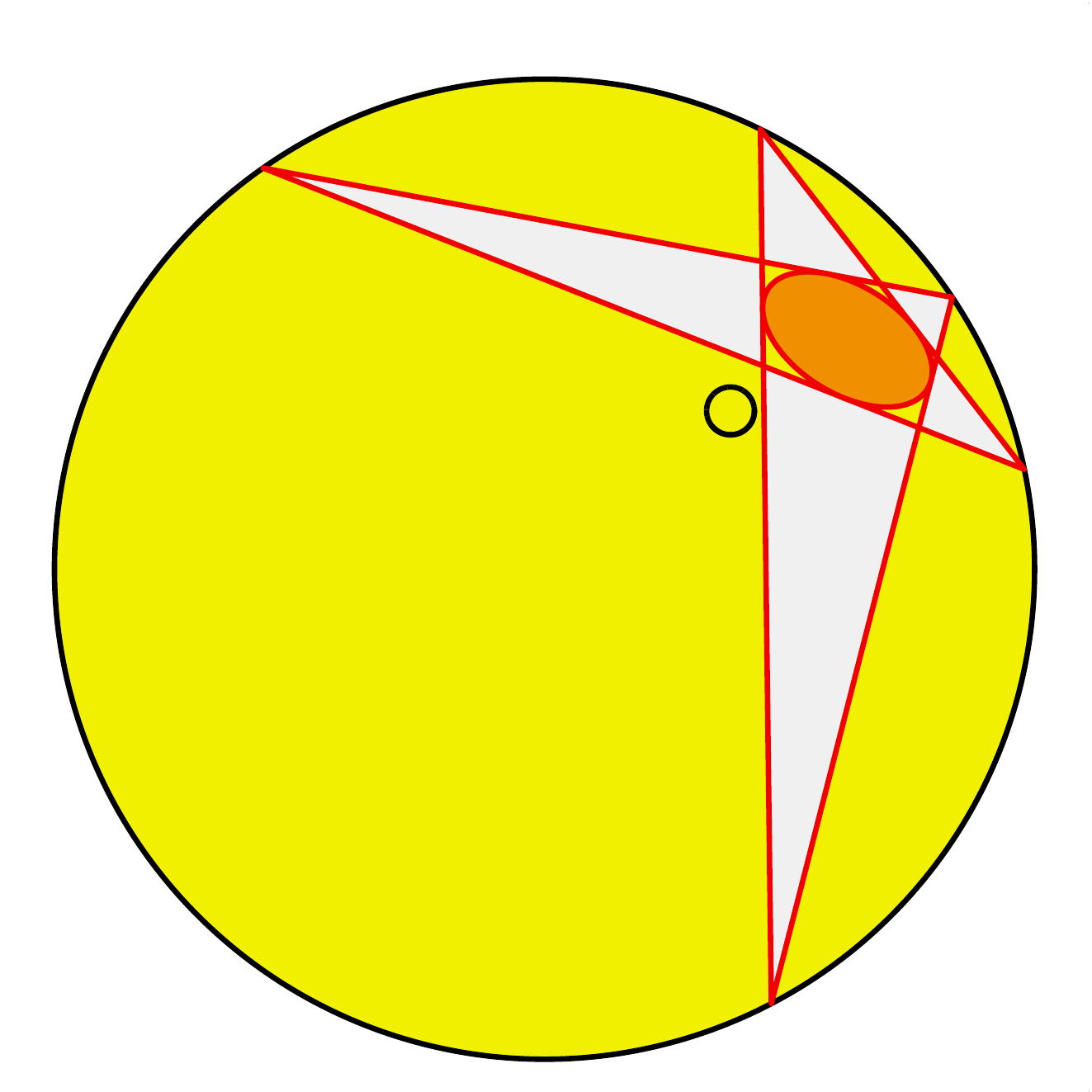}\quad
\includegraphics[height=2.3in]{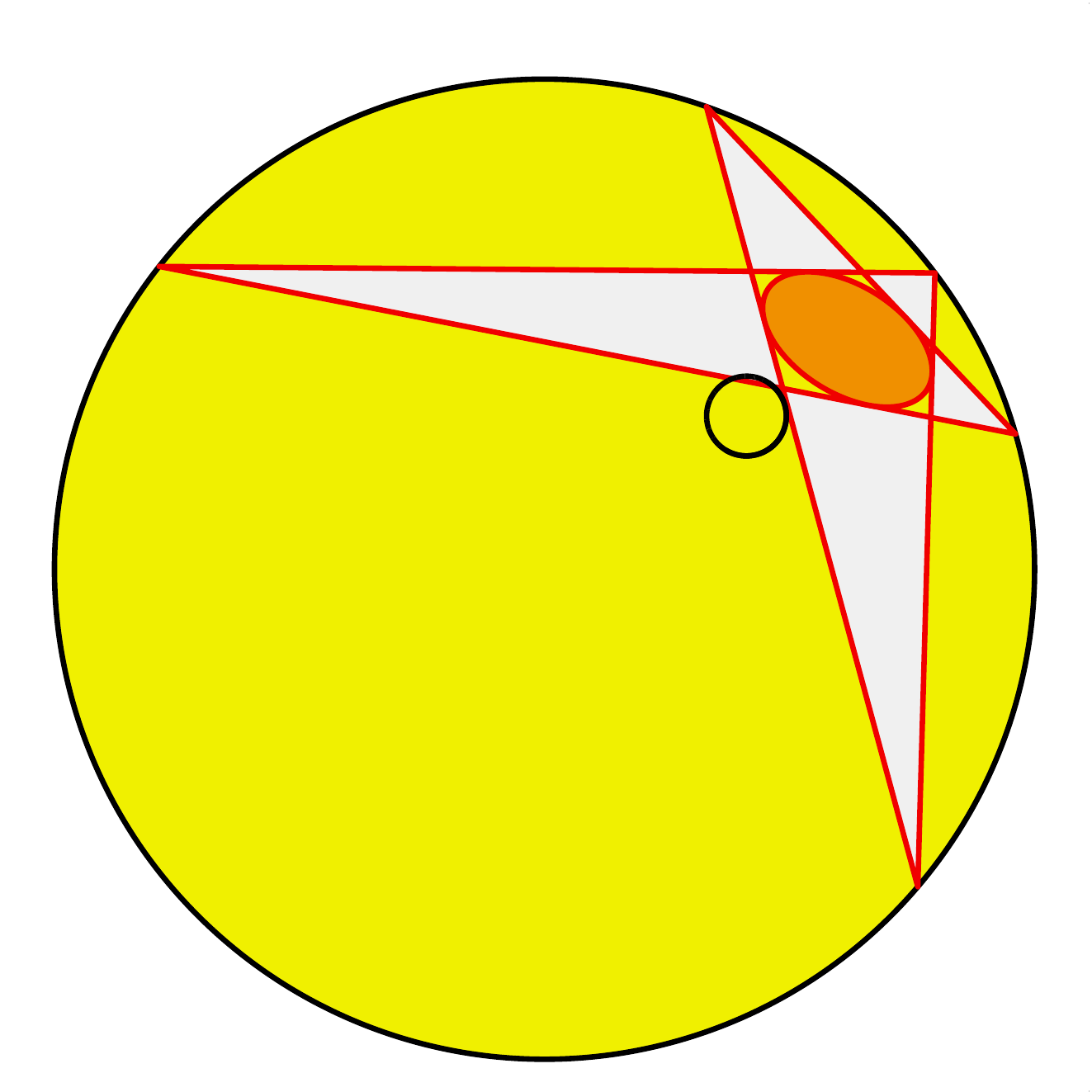}
\caption{$CM_0(P_t)$ and $CM_2(P_t)$.}
\label{experiment}
\end{figure}

We wrote a computer program which tests Theorem \ref{main}.
Figure \ref{experiment} shows two pictures from this program.  The left side
deals with $CM_0$ and the right side deals with $CM_2$.
In both cases, the loci $CM_0(P_t)$ and $CM_2(P_t)$ are black
circles.  The two Poncelet polygons shown belong to the
same Poncelet family, but they are not the same polygon.
The two circles $CM_0(P_t)$ and $CM_2(P_t)$ are
different circles. We are not sure how they are
related to each other.

\section{A  proof} \label{mainpf}

Now we sketch a proof of Theorem \ref{main}. Let us start with the algebraic geometry proof of the Poncelet porism, see \cite{GH}, or \cite{Fl}, for a detailed account.

One starts with complexifying and projectivizing: we assume that $\gamma$ and $\Gamma$ are complex conics in general position in the complex projective plane $\CP^2$ obtained from the affine plane $\CC^2$ by adding a line at infinity.
A complex conic is isomorphic to Riemann sphere, $\CP^1$. 

Let $E$ be the set of flags $(x,L)$ consisting
 of a point $x \in \Gamma$ and a line $L$ through 
$x$ that is tangent to $\gamma$.  The space
$E$ is naturally a Riemann surface, and
the projection $p:E \to \Gamma$ that takes 
$(x,L)$ to point $x$ is a 2-fold branched covering with 
four branch points. These branch  points, invisible in  Figure 5, 
are the four intersection points of $\Gamma$ and $\gamma$. 
One easily computes that the Euler characteristic of 
$E$ is zero. Thus $E$ is a Riemann surface which,
topologically, is a torus.

\begin{figure}[ht]
\centering
\includegraphics[height=1.5in]{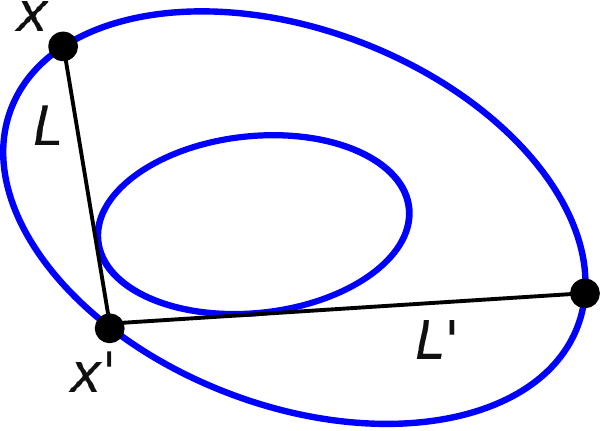}
\caption{Two involutions of the flag space $E$.}
\label{involutions}
\end{figure}

One has two involutions on $E$ 
depicted in Figure \ref{involutions}.
$$
\sigma (x,L) = (x',L),\ \ \tau(x',L) = (x',L').
$$
The Poncelet map $T=\tau \circ \sigma$ is a 
translation of $E$. If a translation of a torus 
has a periodic point of order $n$, then $T^n$ is the 
identity, and this proves the Poncelet theorem.

The two claims of Theorem \ref{main} are proved similarly, 
so we let $CM$ mean either of the two centers of mass involved. 
The group generated by $\sigma$ and $\tau$ is the dihedral group $D_n$. The coordinates $(x,y)$ of the center of mass $CM$ of the Poncelet polygons are $D_n$-invariant meromorphic functions on $E$. Indeed these coordinates are rational functions of the coordinates of the vertices of the polygons (Lemma \ref{CMcoord}), which, in turn, are rational functions of the rational parameter on the conic, which is a meromorphic function $p$ on $E$.

Without loss of generality, assume that $\Gamma$ is a circle. From the complex projective point of view, a circle is a conic that passes through the two circular points $(1:\pm i:0) \in \CP^2$ at infinity. 

When does either coordinate $x$ or $y$ of $CM$ go to infinity? It happens exactly when one of the vertices of the Poncelet $n$-gons   coincides with one of the two circular points. 
This occurs twice for each vertex, with multiplicity one,\footnote{We assume that the conics are in general position.} and then  both $x$ and $y$ go to infinity. Therefore the functions $x$ and $y$ have $4n$ simple poles on $E$, and these poles  comprise two orbits of the group $D_n$, corresponding to the two circular points.

Let $A$ be a pole corresponding to point $(1:i:0)$, and let $z$ be a local holomorphic parameter on $E$ at $A$. Then, at $A$,
$$
x(z) = \frac{a_1}{z}+ a_2 + \ldots,\quad y(z) = \frac{ia_1}{z}+ a_3 + \ldots
$$
where $a_i \in \CC$, and the dots denote  terms of degree one or higher in $z$. Likewise, let $B$ be a pole corresponding to  $(1:-i:0)$. We have, in a local parameter $w$ at $B$:
$$
x(w) = \frac{b_1}{w}+b_2 + \ldots,\quad y(w) = \frac{-ib_1}{w}+ b_3 + \ldots
$$

We want to find $u,v\in\CC$ such that the function $F=(x-u)^2+(y-v)^2$ has no poles at $A$ and $B$. Expanding out $F$ at $A$ and $B$ yields:
$$
F=\frac{2a_1}{z}\left((a_2-u)+i(a_3-v)\right) + \ldots \  {\rm and}\ \  F=\frac{2b_1}{w}\left((b_2-u)-i(b_3-v)\right) + \ldots
$$
Equating both parentheses to zero
 yields two linear equations on $u$ and $v$ which we can
easily solve.  

Consider the function $F$ when
$u$ and $v$ are the solutions to these simultaneous equations.
The function $F$ is $D_n$-invariant and therefore has no poles at all. 
Hence $F$ is holomorphic on $E$.  But the only holomorphic
functions on a compact Riemann surface are constants.  So,
$(X-u)^2+(Y-v)^2=r$, where $r$ is some constant. Since
$CM(P_t)$ contains a curve of points in $\R^2$, the
three constands $u,v,r$ must be real.  In particular,
the complexified locus intersects $\R^2$ in the circle
of radius $r$ centered at $(u,v)$.

\begin{remark} \label{rmk}
{\rm
(i)
Here is a variant of
Theorem \ref{main}.  Consider the
polygons $Q_t$ whose vertices are the tangency 
points of the sides of the Poncelet polygons $P_t$ 
with the ellipse $\gamma$. Then the loci $CM_0(Q_t)$ 
and $CM_2(Q_t)$ are ellipses homothetic to $\gamma$. This is because
the polygon $Q_t$ is also Poncelet:  its sides are 
tangent to the conic which is polar dual to $\Gamma$ 
with respect to $\gamma$.\footnote{We thank A. Akopyan for this argument.} 
\newline
(ii) There is a wealth of results about the loci of various triangle centers of Poncelet triangles. The interested reader is referred to \cite{Ro} for the incenters of 3-periodic billiard trajectories inside an ellipse, to 
\cite{ZKM} for the orthocenters, Gergonne, Nagel, and Lemoine points in bicentric triangles, and to \cite{Sk} for the locus of 
the isogonal conjugate points to a fixed point with respect to the Poncelet triangles inscribed in a circle.
}
\end{remark}

\section{Weill's Theorem} \label{Weill}

Here is Weill's Theorem, stated in the language of Theorem \ref{main}.

\begin{theorem}[Weill]
Let $\gamma \subset \Gamma$ be a pair of nested ellipses that admit a 1-parameter family of Poncelet $n$-gons $P_t$, and let $Q_t$ be the polygons fwhose vertices are the tangency points of the sides of $P_t$ with $\gamma$.  Suppose that both ellipses are homothetic to each other. Then the locus $CM_0(Q_t)$ is a single point.
\end{theorem} 

We present a proof adapted from \cite{MC}. 
We start with Bertrand's proof of the Poncelet 
theorem in the case when both conics are circles (see, e.g., \cite{Fl} or \cite{Sh}).)
There is a subtlety here that we want to emphasize.
For the proof of Poncelet's theorem, which belongs entirely to the domain
of projective geometry, one can always
normalize that both conics are circles.  However, for Theorem \ref{main}
one cannot do this because the various centers of mass are not 
projectively invariant.  

\begin{figure}[ht]
\centering
\includegraphics[height=3in]{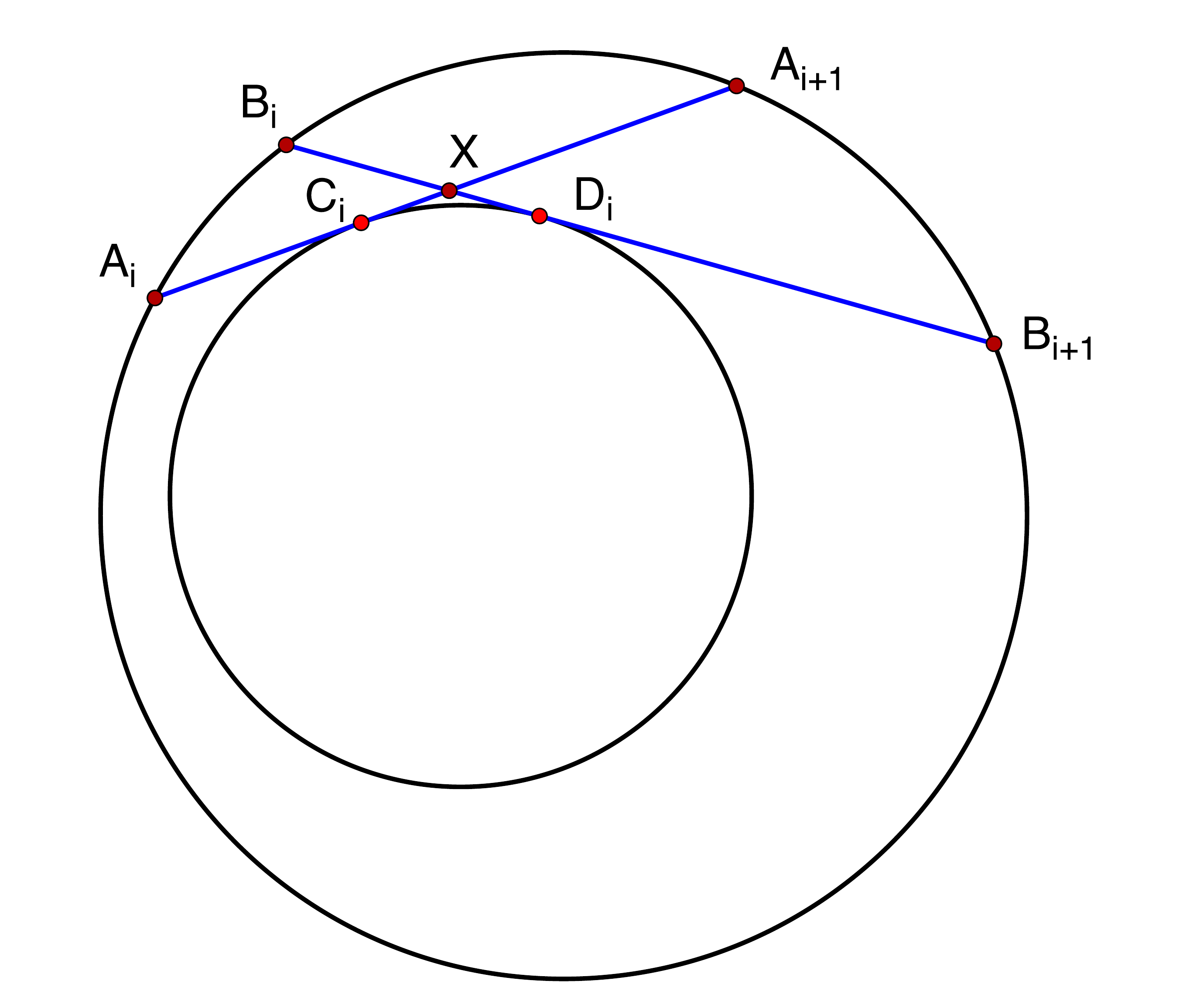}
\caption{Invariant measure for the Poncelet map.}
\label{Bertrand}
\end{figure}

Consider Figure \ref{Bertrand}. Call the map $A_i \mapsto A_{i+1}$ the Poncelet map. The proof consists of constructing a cyclic coordinate $t$ on the outer circle such that, in this coordinate, the Poncelet map is $t \mapsto t+c$, where the constant $c$ depends on the mutual position of the circles. The periodicity property of such a map depends only on $c$, and not on the starting point. 

What follows is, essentially, the argument 
from Theorem XXX, figure 102, in I. Newton's ``Principia" where Newton studies the gravitational attraction of spherical bodies. 

Assume that points $A_i$ and $B_i$ are infinitesimally close (we keep them apart to avoid cluttering the figure). The triangles $A_iB_iX$ and $A_{i+1} B_{i+1}X$ are similar, and hence 
$$
\frac{|A_i B_i| }{|A_iX|} = \frac{|A_{i+1} B_{i+1}|}{|B_{i+1}X|}. 
$$
Denote by $F(A)$ the length of the tangent segment from a point $A$ on the outer circle $\Gamma$ to the inner circle $\gamma$. Then, in the limit $B_i \to A_i$, 
$$
\frac{dA_i}{F(A_i)} = \frac{dA_{i+1}}{F(A_{i+1})},
$$
where $d A$ is the angular measure on the circle. It follows that  the measure $dt(A) := dA/F(A)$ is invariant under the Poncelet map. Hence the map is a  translation in the $t$-coordinate on $\Gamma$.

Now we are ready to prove Weill's Theorem. Denote the infinitesimal angle $A_iXB_i$ by $\varepsilon$. Up to infinitesimals, the angles made by the segments $A_i A_{i+1}$ and $B_i B_{i+1}$ with the circle are equal; let $\varphi$ denote this angle. Let the outer circle be unit, and the inner have radius $R$. Then, in the limit $B_i \to A_i$,
$$
\frac{|C_iD_i|}{|A_iA_{i+1}|} = \frac{\varepsilon R}{2\sin\varphi} = \frac{R}{2} \frac{|A_iB_i|}{|A_iX|} = \frac{R}{2}\, dt(A_i),
$$
where the second equality follows from the Sine Rule for triangle $A_iXB_i$.

Consider an infinitesimal motion of the Poncelet polygon $P=(\ldots A_i A_{i+1}\ldots)$ to $P'=(\ldots B_i B_{i+1} \ldots)$. The center of mass of the points $C_i$ is displaced by the infinitesimal  vector
$$
\sum_i C_i D_i = \frac{R}{2}\, dt\, \sum_i A_i A_{i+1},
$$
and the latter sum is zero because the polygon $P$ is closed.
Therefore the velocity of the center of mass vanishes, and it remains stationary.
\medskip

\noindent
{\bf Acknowledgments}. We are deeply grateful to Mr. X for sharing with us the historical letter from his family archive. Many thanks to A. Akopyan, I. Dolgachev, D. Ivanov, J. Jeronimo, S. Lvovsky, S. Markelov for useful discussions, and  to A. Zaslavsky for making available their harder-to-find article \cite{ZC}. 
The authors were inspired by Martin Gardner's essay \cite{Ga}.
The first author was supported by NSF grant DMS-1204471.
The second author was supported by NSF grants DMS-1105442 and DMS-1510055.

\end{document}